\documentclass[10pt]{article}

\newcommand{\R}{{\sf R\hspace*{-0.9ex}\rule{0.15ex}%
{1.5ex}\hspace*{0.9ex}}}

\newcommand{\N}{{\sf N\hspace*{-1.0ex}\rule{0.15ex}%
{1.3ex}\hspace*{1.0ex}}}

\newtheorem{theorem}{Theorem}[section]

\newtheorem{definition}{Definition}

\begin{document}

\title{Stochastic extrema as stationary phases of characteristic functions}

\author{Sergey Nikitin \\
\thanks{Department of Mathematics and Statistics, Arizona State University, Tempe, AZ 85287-1804  {\tt \small nikitin@asu.edu}}
 }

\maketitle

\begin{abstract}

The paper is dealing with semi-classical asymptotics of a characteristic function for a stochastic process. The main technical tool is provided by the stationary phase method. The extremal range for a stochastic process is defined by limit values of the complex logarithm of the characteristic function. The paper also outlines a numerical method for calculating  stochastic extrema.  

\end{abstract}

\section{INTRODUCTION}

The extremum for a stochastic process admits transparent numerical presentation in terms of limit set of its characteristic function which is treated as a high frequency integral \cite{Guillemin}, \cite{McClure}, \cite{Maslov}, \cite{Maslov_Fedoriuk}. The proposed concept of stochastic extremum is compatible with other known methods of assessing extrema of stochastic functions (see, e.g., \cite{Bohachevsky}, \cite{Brooks}, \cite{Burry}, \cite{Drees}, \cite{Gumbel}, \cite{PeterHall}). However, the ideology of high frequency integrals, though different, is close to the simulated annealing technique \cite{Bohachevsky}, \cite{Brooks}. In our approach the role of the parameter like the inverse of "temperature" (from the annealing process ) is played by the frequency and in order to calculate the extremum we increase the frequency. Method of high frequency integrals is able to calculate extrema of stochastic processes with several variables, and therefore can be applied to analyze images and three-dimensional data samples. When it is possible to apply the method of assessing extrema with principal component functions (in the sense of Karhunen-Lo\`eve representation) \cite{PeterHall}, then our approach works as well and practically leads to the same results. However, it does not rely on any type of Karhunen-Lo\`eve representations and in this sense our concept of stochastic extremum is of equal or more general nature than the Karhunen-Lo\`eve representation itself. Moreover, the method of treating stochastic extrema as stationary phases of its characteristic function leads us to a transparent numerical procedure that allows efficient estimation of the stochastic extremum and evaluation of its statistical significance.

\section{STATIONARY PHASE}

\vspace{0.1cm}

Our goal is to introduce the concept of extremum for a stochastic process. In order to do that we employ the high-frequency integrals:
$$
I(k,\omega) = \int_{-\infty}^{\infty} \varphi (t,\omega) e^{ik\cdot f (t,\omega)} dt,
$$
where $\omega \in \Omega$ is a parameter; both $f$ and $\varphi$ are real functions that are infinitely many times differentiable with respect to $t.$ Moreover, $\varphi (t,\omega)$ has a finite time support for any fixed $\omega \in \Omega, \;\;supp_t (\varphi)$ is a subset of a closed interval from $\R.$ Throughout the paper $\R$ denotes the set of real numbers. For reader convenience, we recall the basic properties of the high-frequency integrals (for further reading on this subject see, e.g. \cite{McClure}, \cite{Maslov}, \cite{Maslov_Fedoriuk}).

\vspace{0.1cm}

If $supp_t (\varphi) \subset [ a, b]$ and 
$$
\frac{d}{dt} f(t,\omega) \not= 0\;\;\forall\;t \in [ a, b]\;\;\mbox{ and } \;\;\forall\; \omega \in \Omega
$$
then
$$
I(k,\omega) = \int_{-\infty}^{\infty} \varphi (t,\omega) e^{ik\cdot f (t,\omega)} dt=\int_{a}^{b} \varphi (t,\omega) e^{ik\cdot f (t,\omega)}dt 
$$
and integrating by parts $n$ times yields
$$
I(k,\omega) =(\frac{1}{ik})^n\int_{a}^{b} L^n(\varphi) (t,\omega) e^{ik\cdot f (t,\omega)}dt 
$$
where the linear operator $L$ is defined as
$$
L(\varphi) = -\frac{d}{dt}(\frac{\varphi}{f_t} )
$$
and $f_t$ denotes the derivative of $f(t,\omega)$ with respect to time,
$$
f_t (t,\omega)= \frac{d}{dt} f(t,\omega).
$$
As one can see, if the phase $f(t,\omega)$ does not have critical points in $supp_t(\varphi)$ then
$$
I(k,\omega) = O(\frac{1}{k^n}) \;\;\forall\;n \in \N \;\;\mbox{ and }\;\;\forall \omega \in \Omega,
$$
where $\N$ denotes the set of natural numbers.  This fact is often represented as
$$
I(k,\omega) = O(\frac{1}{k^\infty})  \;\;\forall \omega \in \Omega.
$$ 
The critical (or stationary) points of the phase $f(t,\omega)$ make the main contribution into the high-frequency integral $I(k,\omega)$ as $k\to\infty .$
\begin{definition}
A point $(t^\star, \omega^\star) \in \R \times \Omega$ is called a stationary phase (point) if 
$$
\frac{d}{dt}f(t^\star , \omega^\star ) = 0.
$$
\end{definition}

The set of all stationary phase points for $f$ is addressed as $St(f) \subset \R \times \Omega.$ Now let us turn our attention to calculating the contribution of a stationary phase point $(t^\star, \omega^\star) \in St(f)$ into the high-frequency integral $I(k,\omega).$ A stationary phase point $(t^\star, \omega^\star) \in St(f)$ is said to have order $m\in \N$ if $m$ is the first natural number for which
$$
(\frac{d}{dt})^m f(t^\star, \omega^\star ) \not= 0.
$$ 
The set of such stationary phase points is denoted by $St^m(f).$
The Taylor expansion near $t^\star$ is
$$
f(t,\omega^\star ) - f(t^\star,\omega^\star ) = \frac{1}{m!}(\frac{d}{dt})^m f(t^\star ,\omega^\star) (t-t^\star)^m + O( (t-t^\star)^{m+1} )  \mbox{ as } \;t\to t^\star.
$$
Consider the change of coordinates
$$
x (t,\omega^\star) = (sign(f^{(m)}_t(t^\star ,\omega^\star)  ) \cdot (f(t,\omega^\star ) - f(t^\star,\omega^\star )  ) )^{\frac{1}{m}},
$$
where $sign(f^{(m)}_t(t^\star ,\omega^\star))$ denotes the sign of
$$
(\frac{d}{dt})^m f(t^\star ,\omega^\star) .
$$
Since 
$$
\frac{d}{dt} x (t,\omega^\star) = \vert \frac{1}{m!}(\frac{d}{dt})^m f(t^\star ,\omega^\star) \vert ^{\frac{1}{m}} + O(t-t^\star)
$$
the change of coordinates is not degenerate on some interval $Q_\varepsilon ,$
$$
t^\star - \varepsilon < t < t^\star + \varepsilon 
$$
Let us take infinitely many times differentiable function $h$ with  $supp_t(h) \subset Q_\varepsilon$ (the set of such functions is denoted as $C_0^\infty (Q_\varepsilon)).$ Assume also that $h( t, \omega^\star)=1$ in a neighborhood of $t^\star.$ Then 
$$
I(k,\omega^\star ) =  \int_{-\infty}^{\infty} \varphi (t,\omega^\star ) h(t,\omega^\star )   e^{ik\cdot f (t,\omega^\star)} dt + \int_{-\infty }^{\infty } \varphi (t,\omega^\star ) (1 -  h(t,\omega^\star )) e^{ik\cdot f (t,\omega^\star)}dt 
$$
and in order to find the contribution of the stationary phase $ (t^\star , \omega^\star )$ we need to calculate asymptotics for
$$
\int_{-\infty}^{\infty} \varphi (t,\omega^\star ) h(t,\omega^\star )   e^{ik\cdot f (t,\omega^\star)} dt.
$$
After making the change of coordinates $x= x (t,\omega^\star)$ in the integral 
$$
I(k,\omega^\star )= e^{ikf(t^\star,\omega^\star)} \cdot \int_{-\infty}^{\infty} \varphi (t,\omega^\star ) h(t,\omega^\star) e^{ik(f(t,\omega^\star) -f(t^\star,\omega^\star)) }dt
$$
we have
$$
I(k,\omega^\star ) = e^{ikf(t^\star,\omega^\star)} \cdot \int_{-\infty}^{\infty}  \varphi (x,\omega^\star ) \frac{h(x,\omega^\star) }{x_t}e^{sign(f^{(m)}_t(t^\star ,\omega^\star)) \cdot  ikx^m}dx,
$$
where $x_t$ denotes $\frac{d}{dt} x(t,\omega^\star).$ The integral 
$$
\int_{-\infty}^{\infty} e^{\pm ikx^m}dx
$$
can be calculated by reducing it to the linear combination of the integrals like$$
\int_{0}^{\infty} e^{\pm ikx^m}dx
$$
and then evaluating the latter with the help of the integral along the curve on a complex plane  \cite{Maslov_Fedoriuk}. The curve consists out of the segment of $x$-axis $0\le x \le \rho,$ the arc of a circle 
$$
\rho \cdot e^{\pm i\tau}\;\;(0\le \tau \le \frac{\pi}{2m})
$$
 and the segment of the straight line 
$$
r\cdot e^{\pm i\frac{\pi}{2m}}\;\;(\rho \ge r \ge 0).
$$
 Taking $\rho \to \infty$  yields that 
$$
\int_{-\infty}^{\infty} e^{\pm ikx^m}dx = \cos(\frac{\pi}{2m})\frac{1}{k^\frac{1}{m}} \cdot C_m \;\;\mbox{ for odd } m
$$
 and
$$
\int_{-\infty}^{\infty} e^{\pm ikx^m}dx = \frac{e^{\pm i\frac{\pi}{2m} }}{k^\frac{1}{m}} \cdot C_m\;\;\mbox{ for even } m
$$
where
$$
C_m = 2\cdot \int_{0}^\infty e^{-x^m} dx. 
$$
Taking into account that $h(t^\star,\omega^\star)=1,$ for even $m$ we have
$$
I(k,\omega^\star ) = \varphi(t^\star,\omega^\star) \cdot C_m \cdot \big(\frac{m!}{k \cdot \vert f^{(m)}_t(t^\star ,\omega^\star) \vert} \big) ^{\frac{1}{m}}\cdot e^{sign(f^{(m)}_t(t^\star ,\omega^\star)) i\frac{\pi}{2m} }    \cdot e^{ikf(t^\star,\omega^\star)} + 
$$
$$
  e^{ikf(t^\star,\omega^\star)}  \int_{-\infty}^{\infty} (\varphi (x,\omega^\star ) \cdot  \frac{h(x,\omega^\star) }{x_t} - \big( \varphi (x,\omega^\star )   \frac{h(x,\omega^\star) }{x_t}     \big)\big\vert_{x=0} )  e^{sign(f^{(m)}_t(t^\star ,\omega^\star)) \cdot  ikx^m}dx 
$$
and the latter integral has the asymptotic
$$
O(\frac{1}{k^\frac{2}{m}})  \;\;\mbox{ as }\;\;k\;\to \; \infty
$$
 as long as there are no other stationary phase points present.
If $m$ is odd then
$$
I(k,\omega^\star ) = \varphi(t^\star,\omega^\star) \cdot C_m \cdot \big(\frac{m!}{k \cdot \vert f^{(m)}_t(t^\star ,\omega^\star) \vert} \big) ^{\frac{1}{m}}\cdot \cos(\frac{\pi}{2m} )    \cdot e^{ikf(t^\star,\omega^\star)} + O(\frac{1}{k^\frac{2}{m}})
$$
as $\;\;k\;\to \; \infty.$

\vspace{0.1cm}

In conclusion to this section we formulate the basic results of stationary phase method in the form of the following formal statement.
\begin{theorem}
\label{stateionaryPhase}
Let $\varphi (t, \omega),\;\;f(t, \omega)$ be real functions that are infinitely many times differentiable with respect to $t$ for any $\omega \in \Omega.$ Moreover, for any  $\omega \in \Omega$ one can find an interval $[a(\omega),b(\omega)]\subset \R$ such that 
$$
\varphi (t, \omega) \in C^\infty_0([a(\omega),b(\omega)]).
$$
Then the following statements hold.
\begin{itemize}
\item[i.] If $[a(\omega),b(\omega)] \cap St(f) = \emptyset$ then
$$
I(k,\omega) = O(\frac{1}{k^\infty})\;\;\mbox{ as } \;\;k\to\infty
$$
\item[ii.]  If $[a(\omega),b(\omega)] \cap St(f) = \{t_j\}$ then $I(k,\omega)$ has the asymptotic
$$
\sum_{\mbox{even } m_j }\left(\varphi(t_j,\omega) \cdot C_{m_j} \cdot \big(\frac{m_j!}{k \cdot \vert f^{(m_j)}_t(t_j ,\omega) \vert} \big) ^{\frac{1}{m_j}}\cdot e^{sign(f^{(m_j)}_t(t_j ,\omega)) i\frac{\pi}{2m_j} }    \cdot e^{ikf(t_j,\omega)} \right. +
$$
$$
\left. O(\frac{1}{k^\frac{2}{m_j}})\right) +  \sum_{\mbox{odd } m_j } \left(\varphi(t_j,\omega) \cdot C_{m_j} \cdot \big(\frac{m_j!}{k \cdot \vert f^{(m_j)}_t(t_j ,\omega) \vert} \big) ^{\frac{1}{m_j}}\cdot \cos(\frac{\pi}{2m_j} )    \cdot e^{ikf(t_j,\omega)} \right. +
$$
$$
\left. O(\frac{1}{k^\frac{2}{m_j}}) \right)
$$
as $k\to \infty .$
\end{itemize}
\end{theorem}

\section{STOCHASTIC EXTREMUM}
Consider a real valued stochastic process $\xi(t) $ defined in the probability space $(\Omega, A(\Omega), P),$ where $A(\Omega)$ is a $\sigma$-algebra of subsets from $\Omega$ and $P$ is a measure of probability. Throughout the paper we assume that $\xi(t)$ takes only non-negative real values. In the context of this paper it is tacitly assumed that  
$$
     \xi(t, \omega) = g(t,\omega),
$$
where $\omega$ is a stochastic variable which does not depend on time and $g(t,x)$ is a smooth non-negative function of its arguments.

\vspace{0.1cm}

In order to examine whether $\xi(t) $ has extrema on interval $[a,b]\subset \R$ we take a function $\varphi_\varepsilon (t)\in C^\infty_0([a-\varepsilon , b+\varepsilon])$ such that $\varepsilon > 0,\;\;\;\varphi_\varepsilon (t)\ge 0\;\;\forall \;t\in \R$ and
$$
  \varphi_\varepsilon (t) =1 \;\;\forall\; t\in [a,b].
$$
Then we analyze the asymptotics of the high-frequency integral
$$
 \int_\Omega \int_{-\infty }^{\infty} \varphi_\varepsilon (t) e^{ik \xi(t,\omega)}dtdP(\omega) .
$$
The formal description of the situation when $\xi(t)$ does not have stochastic extrema on $[a,b]$ sounds as follows.

\begin{definition}
\label{noExtrema}
Let $\ln(z)$ denote a fixed branch of the complex logarithm. Then a real-valued stochastic process $\xi(t)$ does not have stochastic extrema on $[a,b]$ if one can find a positive real number $\varepsilon$ such that 
$$
\lim_{k\to \infty } Re\{\frac{1}{ik} \ln( \int_\Omega \int_{-\infty }^{\infty} \varphi_\varepsilon (t) e^{ik \xi(t, \omega )}dt dP(\omega)  )\}= 0.
$$ 
\end{definition}
The logical negation of this statement describes the intervals where stochastic extrema for $\xi(t)$ occur. In other words, the interval  $[a,b]$ contains stochastic extrema if $\forall\;\varepsilon >0$
$$
\lim_{k\to \infty }Re\{\frac{1}{ik} \ln(\int_\Omega  \int_{-\infty }^{\infty} \varphi_\varepsilon (t) e^{ik \xi(t,\omega )}dt dP(\omega) )\}\not= 0.
$$
Sometimes it is possible to estimate extremal values for $\xi(t)$ on $[a,b]$ with the help of the following 
\begin{equation}
\label{Smax}
Smax_\varepsilon ([a,b]) =  \overline{\lim}_{k\to \infty }Re \left\{\frac{1}{ik} \ln( \int_\Omega \int_{-\infty }^{\infty} \varphi_\varepsilon (t) e^{ik \xi(t,\omega)}dt dP(\omega)  )\right\}
\end{equation}
and
\begin{equation}
\label{Smin}
Smin_\varepsilon ([a,b]) =  {\underline \lim}_{k\to \infty }Re\left\{\frac{1}{ik} \ln(\int_\Omega   \int_{-\infty }^{\infty} \varphi_\varepsilon (t) e^{ik \xi(t,\omega)}dt dP(\omega) )\right\},
\end{equation}
where $\overline{\lim}$ and ${\underline \lim}$ denote upper and lower limits, respectively. 

\vspace{0.1cm}

For majority of applications $Smax_\varepsilon ([a,b]),\;\;Smin_\varepsilon ([a,b])$ deliver boundaries for extremal values of a real stochastic process $\xi(t).$ One can justify that when the stochastic process $\xi(t)$ has additional properties like, for example, $\Omega$ is a smooth manifold and $\xi(t,\omega )$ is a smooth function of $t$ and $\omega.$ Let  $\frac{dP}{d\omega}$ denote the probability density function for $P(\omega).$ Then the following statement takes place.

\begin{theorem}
\label{smoothExpectation}
Let $\Omega$ be a smooth manifold. Assume also that density  $\frac{dP}{d\omega}$ exists; $\frac{dP}{d\omega} \in C^\infty (\Omega),$ $\xi(t,\omega )\in C^\infty(\R ,\Omega )$ and there is only one stationary phase point $( t^\star, \omega^\star) \in St^2(\xi)\cap \left( [a-\varepsilon ,b + \varepsilon ]\times \Omega \right)$ with $\frac{\partial}{\partial \omega} \xi(t^\star, \omega^\star ) = 0 $ and  the second derivative $\frac{\partial^2 \xi }{\partial t \partial \omega}$ has non-zero determinant
\begin{equation}
\label{nonZero2nd}
J(t^\star, \omega^\star)=\det \left( \frac{\partial^2 \xi }{\partial t \partial \omega} (t^\star, \omega^\star) \right) \not= 0
\end{equation}
at  the stationary phase point $(t^\star, \omega^\star).$
If $\frac{d}{d\omega}P(\omega^\star) >0$ then the following is true:
$$
Re\left\{ \frac{1}{ik} \cdot \ln \left( \int_\Omega  \int_{-\infty }^{\infty} \varphi_\varepsilon (t) e^{ik \xi(t,\omega)}dt  dP(\omega) \right) \right\} = \xi(t^\star,\omega^\star ) + O(\frac{1}{k})\;\;\;\mbox{ as }\;\;k\to \infty.
$$
\end{theorem}

\vspace{0.1cm}

{\bf Proof.}

\vspace{0.1cm}

Let $n$ denote the dimension of $\Omega.$ Then applying the stationary phase method \cite{Maslov} ,\cite{Maslov_Fedoriuk} to 
$$
 \int_\Omega \int_{-\infty }^{\infty} \varphi_\varepsilon (t) e^{ik \xi(t,\omega )}dt dP(\omega)
$$
yields
$$
\int_\Omega \int_{-\infty }^{\infty} \varphi_\varepsilon (t) e^{ik \xi(t)}dt dP(\omega) =  
$$
$$
\frac{d}{d\omega}P(\omega^\star) \cdot (\frac{2\pi}{k})^{\frac{n+1}{2}} \cdot \frac{1}{\sqrt{\vert J(t^\star,\omega^\star) \vert}}  \cdot e^{sign(\frac{\partial^2 \xi }{\partial t \partial \omega}(t^\star,\omega^\star )) i\frac{\pi}{4} }    \cdot e^{ik\xi(t^\star,\omega^\star )} +   O(\frac{1}{k^{1+\frac{n+1}{2}}}) 
$$
where 
$$
sign(\frac{\partial^2 \xi }{\partial t \partial \omega}(t^\star,\omega^\star ))
$$
denotes the difference between the number of positive  and the number of negative  eigenvalues of the corresponding quadratic form.
Taking the complex logarithm of 
$$
\frac{d}{d\omega}P(\omega^\star)  (\frac{2\pi}{k})^{\frac{n+1}{2}} \cdot \frac{1}{\sqrt{\vert J(t^\star,\omega^\star) \vert}}  \cdot e^{sign(\frac{\partial^2 \xi }{\partial t \partial \omega}(t^\star,\omega^\star )) i\frac{\pi}{4} }    \cdot e^{ik\xi(t^\star,\omega^\star )} \cdot(1 + O(\frac{1}{k}))
$$
we obtain
$$
\ln (\int_{-\infty }^{\infty} \varphi_\varepsilon (t) e^{ik \xi(t)}dt) = ik\xi(t^\star,\omega^\star) - \frac{n+1}{2}\ln (k) + \ln (\frac{d}{d\omega}P(\omega^\star) ) + \frac{n+1}{2} \cdot  \ln (2\pi) -
$$
$$
\frac{1}{2} \ln\big( \vert J(t^\star,\omega^\star) \vert \big)  + i\frac{\pi}{4} \cdot sign(\frac{\partial^2 \xi }{\partial t \partial \omega}(t^\star,\omega^\star ))+   O(\frac{1}{k})
$$
as $\;\;k\to \infty.$ Dividing by  $ik$ and taking real part complete the proof.
{\bf Q.E.D.}

\vspace{0.1cm}

Formulas (\ref{Smax}), (\ref{Smin}) together with Theorem \ref{smoothExpectation} give us a recipe for assessing extrema of a real stochastic process $\xi(t)$ on an arbitrary time interval $[a,b].$ In order to do that one needs to analyze the limit properties of the following integral for big values of $k:$
$$
Re\left\{\frac{1}{ik}  \int_{\gamma_k} \frac{dz}{z} \right\},
$$
where the complex curve $\gamma_k$ is defined as
$$
\gamma_k = \left\{ \int_\Omega \int_{-\infty }^{\infty} \varphi_\varepsilon (t) e^{i \lambda \xi(t,\omega )}dt dP(\omega);\;\; 0 \le \lambda  \le k  \right\}.
$$
Notice that $\gamma_k$ can be interpreted as the characteristic function of $ \xi(t,\omega ),$ where $\varphi_\varepsilon (t)$ is chosen so that 
$$
\int_{-\infty }^{\infty} \varphi_\varepsilon (t) dt =1.
$$
 In conclusion we notice that the condition (\ref{nonZero2nd}) can be relaxed with the help of Theorem \ref{stateionaryPhase} (see also \cite{Atiyah}, \cite{Malgrange}). In some applications (e.g., when estimating tails of distributions ) it is beneficial that the truncation function $\varphi_\varepsilon (t)$ also depends on $\omega.$

\end{document}